# The Influence of Nitsche Stabilization on Geometric Multigrid for the Finite Cell Method[*]

S. Saberi [†]    G. Meschke [‡]    A. Vogel [†]

**Abstract**

Immersed finite element methods have been developed as a means to circumvent the costly mesh generation required in conventional finite element analysis. However, the numerical ill-conditioning of the resultant linear system of equations in such methods poses a challenge for iterative solvers. In this work, we focus on the finite cell method (FCM) with adaptive quadrature, adaptive mesh refinement (AMR) and Nitsche's method for the weak imposition of boundary conditions. An adaptive geometric multigrid solver is employed for the discretized problem. We study the influence of the mesh-dependent stabilization parameter in Nitsche's method on the performance of the geometric multigrid solver and its implications for the multilevel setup in general. A global and a local estimate based on generalized eigenvalue problems are used to choose the stabilization parameter. We find that the convergence rate of the solver is significantly affected by the stabilization parameter, the choice of the estimate and how the stabilization parameter is handled in multilevel configurations. The local estimate, computed on each grid, is found to be a robust method and leads to rapid convergence of the geometric multigrid solver.

# 1 Introduction

Numerical approximation of partial differential equations on complex domains is a labor-intensive process, where the generation of an appropriate boundary-conforming computational mesh can constitute a salient portion of the simulation workflow. This has motivated the development of a class of advanced discretization methods that can broadly be categorized under the umbrella of immersed or unfitted finite element methods, namely eXtended FEM (XFEM)[1], finite cell method (FCM)[2,3] and cutFEM[4] methods fall under this category. The classification of each method depends on the specific techniques that are used for domain integration and imposition of boundary conditions. Lagrange multipliers,[5–10] the penalty method[11,12], Nitsche's method[13–20], the LS-Nitsche's method,[21,22] etc. are some of the techniques that are commonly used for the imposition of boundary conditions. Both the penalty method and Nitsche's method impose essential boundary conditions in a weak sense, i.e., the variational formulation is modified to account for essential boundary conditions rather than introducing explicit constraints on the state variables. In addition, they introduce no additional unknowns and preserve the symmetry, positive-definiteness and banded structure of the matrix. In contrast to Nitsche's method, the penalty method is not variationally consistent[5,18] and may lead to severe ill-conditioning of the system. The finite

---

[*]Financial support was provided by the German Research Foundation (*Deutsche Forschungsgemeinschaft, DFG*) in the framework of the collaborative research center SFB 837 *Interaction Modeling in Mechanized Tunneling*, grant number 77309832.

[†]High Performance Computing, Ruhr University Bochum, Universitätsstr. 150, 44801 Bochum, Germany

[‡]Institute for Structural Mechanics, Ruhr University Bochum, Universitätsstr. 150, 44801 Bochum, Germany

cell method makes use of adaptive quadrature for domain integration and a weak imposition method for essential boundary conditions,[2,3,23] and Nitsche's method is usually preferred in the context of FCM.[24–26]

From a computational point of view, one of the main hurdles for the employment of such numerical methods is the solution of the discretized problem that has often limited their deployment on large scale problems. While direct solvers are robust and solve the system to machine accuracy, they normally suffer from sub-optimal complexity and become extremely expensive for large systems[27]. On the other hand, iterative solvers can offer much better complexity and concurrency on parallel machines; however, typically their convergence highly depends on the spectrum of the system matrix and is therefore mesh dependent.[27] Multigrid methods have been successfully used to remove the mesh dependence of iterative solvers in the finite element discretization of different classes of PDEs.[28] In the context of unfitted finite element methods, it is well known that small cut fractions can lead to severe ill-conditioning of the system matrix and, therefore, limit the usability of iterative solvers. A geometric multigrid solver for three CutFEM formulations[29] and for XFEM[30] using the Nitsche's method was recently studied. Geometric multigrid has been recently studied for the finite cell method with a penalty formulation[31] and with the Nitsche's method.[32,33]

In this work, we employ a finite cell formulation of the Poisson problem where adaptive quadrature and Nitsche's method are used for domain integration and the imposition of boundary conditions, respectively. We formulate a geometric multigrid solver for the solution of the resultant system and investigate the influence of the stabilization parameter in Nitsche's method on the performance of the geometric multigrid solver for the finite cell formulation. The sensitivity of the solver to variations in the stability parameter is studied and possible methods for the estimation of an appropriate stability parameter are explored. Furthermore, we discuss the implications of the mesh-dependent nature of the stability parameter for multilevel methods.

The rest of this work is organized as follows. The numerical formulation of the model problem, the weak imposition of essential boundary conditions, the finite cell formulation and the geometric multigrid solver are described in detail in Section 2. The geometric multigrid solver is studied using a number of numerical benchmarks in Section 3 and the results are discussed. Finally, conclusions are drawn in Section 4.

## 2 Numerical methodology

### 2.1 Model problem

As the model problem in this work, we use the Poisson equation as a representative elliptic partial differential equation which can be written in strong form as

$$\begin{aligned} -\nabla^2 u &= f & \text{in } \Omega, \\ u &= g & \text{on } \Gamma_D, \\ \nabla u \cdot \boldsymbol{n} &= t & \text{on } \Gamma_N, \end{aligned} \qquad (1)$$

where $u$ is the scalar solution variable, $\Omega$ is the spatial domain, $\Gamma = \Gamma_D \cup \Gamma_N$ is the boundary of the domain, $\Gamma_D$ and $\Gamma_N$ are the Dirichlet and Neumann parts of the boundary, respectively, $f$ is the body force function and $\boldsymbol{n}$ is the unit-length outer normal vector to the boundary. $g$ and $t$ are prescribed functions on $\Gamma_D$ and $\Gamma_N$, respectively. Multiplying the strong form with appropriate test functions, integrating over the domain and using the Green's theorem, the boundary-conforming weak form can be written as:

Find $u \in V$ such that for all $v \in V_0$

$$\int_\Omega \nabla v \cdot \nabla u \, d\boldsymbol{x} - \int_{\Gamma_D} v(\nabla u \cdot \boldsymbol{n}) \, d\boldsymbol{s} = \int_\Omega vf \, d\boldsymbol{x} + \int_{\Gamma_N} vt \, d\boldsymbol{s}, \qquad (2)$$

where $v$ are the test functions and

$$V := \{u \in H^1(\Omega) \mid u = g \text{ on } \Gamma_D\},$$
$$V_0 := \{v \in H^1(\Omega) \mid v = 0 \text{ on } \Gamma_D\}, \tag{3}$$

where $H^1$ is the Sobolov space. We note that the term on $\Gamma_D$ in Equation (2) vanishes in the boundary-conforming case as the essential boundary conditions are included in the Space $V$.

## 2.2 Immersed finite element method

The physical domain $\Omega$ is approximated using a boundary-conforming tessellation in classical finite element methods. However, given that $\Omega$ can be arbitrarily complex, a fact which is directly reflected in the effort required for the generation of a boundary-conforming tessellation of the domain, immersed finite elements method embed $\Omega$ in an embedding domain $\Omega_e$, as shown in Figure 1. The computational domain is thereby extended by the fictitious part $\Omega_e \setminus \Omega$. The embedding of $\Omega$ in $\Omega_e$ necessitates two overarching modifications to the standard finite element weak formulation. On the one hand, the weak formulation must recover the physical domain $\Omega$ which is achieved using a penalization factor as explained in detail below. On the other hand, as opposed to the boundary-conforming case, in which the space $V$ includes the essential boundary conditions, Equation 2 must be modified to impose essential boundary conditions weakly in immersed methods since the physical boundary $\Omega$ is not guaranteed to conform to the boundary of the computational domain $\Omega_e$. We start with the latter, namely the weak imposition of essential boundary conditions, for which we use the symmetric Nitsche's method in this work.

Applying the Nitsche's method, the weak formulation in Equation (2) turns into

$$\int_\Omega \nabla v \cdot \nabla u \, d\boldsymbol{x} - \int_{\Gamma_D} v(\nabla u \cdot \boldsymbol{n}) \, d\boldsymbol{s} - \int_{\Gamma_D} (u-g)(\nabla v \cdot \boldsymbol{n}) \, d\boldsymbol{s} + \int_{\Gamma_D} \lambda v(u-g) \, d\boldsymbol{s} = \int_\Omega vf \, d\boldsymbol{x} + \int_{\Gamma_N} vt \, d\boldsymbol{s}, \tag{4}$$

where $\lambda$ is the scalar stabilization parameter in the Nitsche's method. The third term on the left-hand side of Equation (4) is the symmetric consistency term ensuring that the symmetry of the original weak form is retained. The fourth term on the left-hand side of Equation (4) is the stabilization term ensuring that the boundary conditions are satisfied and that the formulation is stable for a large enough $\lambda$.

The weak formulation in Equation 4 is in the next step extended to the embedding domain as follows:

$$\int_{\Omega_e} \alpha \nabla v \cdot \nabla u \, d\boldsymbol{x} - \int_{\Gamma_D} v(\nabla u \cdot \boldsymbol{n}) \, d\boldsymbol{s} - \int_{\Gamma_D} (u-g)(\nabla v \cdot \boldsymbol{n}) \, d\boldsymbol{s} + \int_{\Gamma_D} \lambda v(u-g) \, d\boldsymbol{s} = \int_{\Omega_e} \alpha vf \, d\boldsymbol{x} + \int_{\Gamma_N} vt \, d\boldsymbol{s}, \tag{5}$$

where

$$\begin{cases} \alpha = 1, & \text{in } \Omega, \\ \alpha = 0, & \text{in } \Omega_e \setminus \Omega, \end{cases} \tag{6}$$

where $\alpha$ is the scalar penalization factor. The physical domain is essentially recovered by penalizing the part of the embedding domain that lies outside of the physical domain. We note that in practice, in the fictitious part $\Omega_e \setminus \Omega$, a small positive, non-zero value $\alpha \ll 1$ is used instead of zero in order to avoid severe numerical ill-conditioning.

Let $\mathcal{T}_h := \{K_i\}_{i=1}^{n_K}$ be a tessellation of the computational domain $\Omega_e$ into a set of $n_K$ compact, connected, Lipschitz sets $K_i$ with non-empty interior and $\mathring{K}_i \cap \mathring{K}_j = \emptyset \,\, \forall \, i \neq j$.

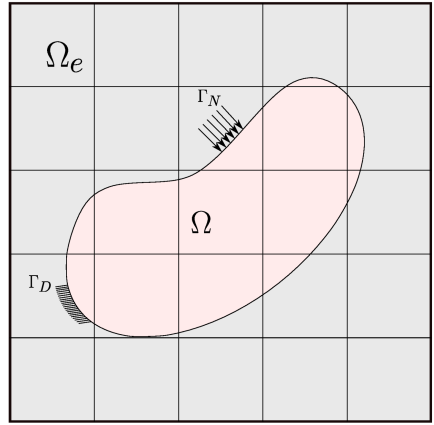

Figure 1: Illustration of the physical domain $\Omega$ along with Dirichlet ($\Gamma_D$) and Neumann ($\Gamma_N$) boundaries in a structured embedding domain $\Omega_e$

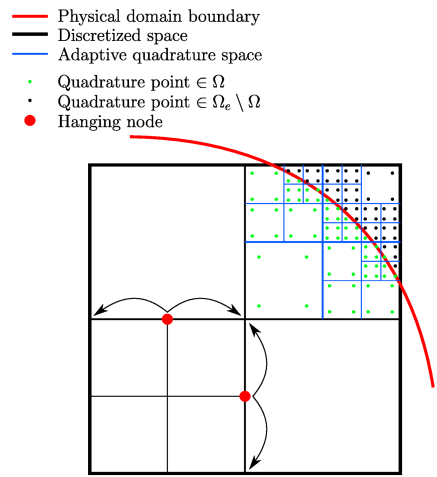

Figure 2: Adaptive quadrature integration on a cell cut by the physical domain and handling of hanging nodes as a result of adaptive mesh refinement

$\overline{\Omega}_{e,h} := \cup_{i=1}^{n_K} K_i$ then defines an approximation of $\Omega_e$. Introducing a finite-dimensional function space $V_h \subset H^1(\Omega_e)$, the following weak formulation is obtained:

Find $u_h \in V_h$ such that for all $v_h \in V_h$

$$a_h(u_h, v_h) = b_h(v_h), \tag{7}$$

with

$$a_h(u_h, v_h) := \int_{\Omega_e} \alpha \nabla v_h \cdot \nabla u_h \, d\boldsymbol{x} - \int_{\Gamma_D} v_h (\nabla u_h \cdot \boldsymbol{n}) \, d\boldsymbol{s} - \int_{\Gamma_D} u_h (\nabla v_h \cdot \boldsymbol{n}) \, d\boldsymbol{s} + \int_{\Gamma_D} \lambda v_h u_h \, d\boldsymbol{s},$$

$$b_h(v_h) := \int_{\Omega_e} \alpha v_h f \, d\boldsymbol{x} + \int_{\Gamma_N} v_h t \, d\boldsymbol{s} - \int_{\Gamma_D} g (\nabla v_h \cdot \boldsymbol{n}) \, d\boldsymbol{s} + \int_{\Gamma_D} \lambda v_h g \, d\boldsymbol{s}. \tag{8}$$

The embedding domain is typically regular and can be generated efficiently using a variety of data structures such as space trees.[34] In this work, we employ adaptive quadrature in conjunction with adaptive mesh refinement (AMR) using space trees for volume integration as illustrated in Figure 2. Adaptive quadrature integration is a natural extension of adaptive mesh refinement on space tree data structures. We note that the adaptive integration space (blue lines in Figure 2) is not a part of the global system of equations. As shown in Figure 2, adaptive mesh refinement leads to non-conforming discretizations, i.e., hanging nodes, that must be properly handled in order to obtain a consistent solution. We constrain hanging nodes to their corresponding non-hanging counterparts and remove them from the global system. We use conventional Lagrange shape functions throughout this work.

The symmetric Nitsche's method is robust and accurate, provided that the stabilization parameter is appropriately chosen; however, an insufficiently-large parameter leads to loss of stability of the discretization, and extremely-large parameters are equivalent to the penalty method, which is stable but can cause severe ill-conditioning in the matrix. The choice of the stabilization parameter is, therefore, crucial to the discretization. A number of methods for the estimation of the stabilization parameter involve the solution of general eigenvalue problems.[15,17,35] Applying Young's inequality with $\varepsilon$, the bilinear form in Equation (8) can be estimated by

$$\begin{aligned} a_h(v_h, v_h) \geq &\int_{\Omega_e} \alpha \nabla v_h \cdot \nabla v_h \, d\boldsymbol{x} - \frac{1}{\varepsilon} \int_{\Gamma_D} v_h v_h \, d\boldsymbol{s} \\ &- \varepsilon \int_{\Gamma_D} (\nabla v_h \cdot \boldsymbol{n})(\nabla v_h \cdot \boldsymbol{n}) \, d\boldsymbol{s} + \int_{\Gamma_D} \lambda v_h v_h \, d\boldsymbol{s}. \end{aligned} \tag{9}$$

Assuming a constant $C$ such that

$$C \int_{\Omega_e} \alpha \nabla v_h \cdot \nabla v_h \, d\boldsymbol{x} \geq \int_{\Gamma_D} (\nabla v_h \cdot \boldsymbol{n})(\nabla v_h \cdot \boldsymbol{n}) \, d\boldsymbol{s}, \tag{10}$$

and assuming that the stabilization parameter $\lambda$ is constant over the integration domain, the following relation is obtained

$$a_h(v_h, v_h) \geq (1 - \varepsilon C) \int_{\Omega_e} \alpha \nabla v_h \cdot \nabla v_h \, d\boldsymbol{x} + (\lambda - \frac{1}{\varepsilon}) \int_{\Gamma_D} v_h v_h \, d\boldsymbol{s}. \tag{11}$$

For the coercivity of the bilinear form, it is required that both $(1 - \varepsilon C)$ and $(\lambda - \frac{1}{\varepsilon})$ be positive. It can then easily be shown that necessarily $\lambda > C$. Under this stability-ensuring constraint, the stabilization parameter $\lambda$ should be chosen as small as possible in order to avoid numerical ill-conditioning. Therefore, in order to find an appropriate value for $\lambda$, a good estimate for the lower bound $C$ has to be obtained. This can be achieved from Equation (10), which allows to compute $C$ via a generalized eigenvalue problem as follows:

$$\boldsymbol{K}\boldsymbol{v} = \Lambda \boldsymbol{M}\boldsymbol{v}, \tag{12}$$

where $\boldsymbol{v}$ and $\Lambda$ are the eigenvectors and eigenvalues, respectively, and using a basis $\{\phi_i\} \subset V_h$, the matrices $\boldsymbol{K}$ and $\boldsymbol{M}$ are formed by

$$\begin{cases} \boldsymbol{K}_{ij} := \int_{\Gamma_D} (\nabla \phi_i \cdot \boldsymbol{n})(\nabla \phi_j \cdot \boldsymbol{n}) \, d\boldsymbol{s}, \\ \boldsymbol{M}_{ij} := \int_{\Omega_e} \alpha \nabla \phi_i \cdot \nabla \phi_j \, d\boldsymbol{x}, \end{cases} \qquad (13)$$

for the model problem. The lower bound $C$ for the stabilization parameter in Nitsche's method can now be chosen as the largest eigenvalue $\max_k \Lambda_k$. Please note that $\boldsymbol{M}$ is integrated over the entire domain in Equation (13). However, the corresponding integral can conservatively be restricted to the part of the domain intersected by the Dirichlet boundary in Equation (10) and consequently in Equation (13).

Given the above reasoning, two approaches for the stabilization parameter estimation can be formulated: The global approach estimates $C$ through the solution of one generalized eigenvalue problem given by Equation (12), where the integration of $\boldsymbol{M}$ is restricted to only those elements $K \in \mathcal{T}_h$ which are intersected by $\Gamma_D$. This provides a single stabilization parameter $\lambda_g$ obtained for the entire domain. Alternatively, the local approach is to restrict the estimate in Equation (10) and analogously the generalized eigenvalue problem in Equation (12) to the domain of each element $K \in \mathcal{T}_h$ cut by the Dirichlet boundary. Thereby, element-wise constant stabilization parameters $\lambda_l^K, K \in \mathcal{T}_h$ with $K \cap \Gamma_D \neq \emptyset$, are computed through the solution of a series of local generalized eigenvalue problems. For shorter notations, we denote the set of element-wise parameters $\{\lambda_l^K\}$ by $\lambda_l$.

In this work, we employ both methods and provide a comparison between the two approaches with respect to their impact on the iterative solver performance in Section 3. From a computational point of view, the main difference is to either solve one large eigenvalue problem or several smaller eigenvalue problems. The solution of the generalized eigenvalue problem (12), however, is not trivial due to the rank deficiency of the matrices. A possible algorithm for solving the problem is singular value decomposition which is employed in this work. Due to the increasing size of the eigenvalue problem for the global approach, it might become prohibitively expensive on very fine meshes, whereas the local approach remains applicable.

## 2.3 Geometric multigrid

The discretized form of the model problem as described in Section 2.2 leads to the following system of equations:

$$\boldsymbol{A}\boldsymbol{x} = \boldsymbol{b}, \qquad (14)$$

where $\boldsymbol{A}$ and $\boldsymbol{b}$ are defined according to the bilinear and linear forms in Equation (8), respectively, and $\boldsymbol{x}$ is the solution vector.

It is well known that the existence of small cut fractions, where the intersection between the physical domain and elements of the computational embedding domain is small, can lead to numerical ill-conditioning in immersed finite element methods, and the condition number of the system matrix can be arbitrarily large; therefore, special treatment is necessary for the successful application of iterative methods. In this work, we use a Schwarz-type smoother for the treatment of cut cells.[31,32] The fundamental idea of the multigrid method, namely smoothing the oscillatory frequencies of the error on the fine grid and approximating the smooth frequencies on the coarse grid remains intact.

### 2.3.1 Components

In its most fundamental form, geometric multigrid requires the following components: a hierarchy of grids $\tau_0, \ldots \tau_n$, where $\tau_n$ is the original fine mesh and $\tau_0$ is the coarsest mesh, transfer operators $\boldsymbol{R}_l$ and $\boldsymbol{P}_l = \boldsymbol{R}_{l+1}^T$ that restrict a vector from level $l$ to $l-1$ and prolongate it from $l$ to $l+1$, respectively, a smoother operator $\boldsymbol{S}$ that is employed in a number of pre- and post-smoothing steps for each level except the coarsest level $\tau_0$ and a base solver that is employed on the coarsest level $\tau_0$.

**Grid Hierarchy** We highlight a few important aspects of the adaptive mesh refinement employed in this work. An exemplary grid refinement is shown in Figure 3 with a four-level grid hierarchy. A 2:1 balance is enforced on all grids, which means that no neighbor quadrants will have refinement levels that are more than two levels apart. Starting from a fine grid $\tau_l$, in order to arrive at the immediate coarse grid $\tau_{l-1}$, a loop through the data structure is performed and all quadrants with the maximum refinement level are coarsened. Applying this algorithm recursively leads to a desired number of nested spaces.

**Transfer operators** Given a vector $\boldsymbol{v}_l$ on level $l$, the restricted vector $\boldsymbol{v}_{l-1}$ on the coarse level $l-1$ can be obtained as
$$\boldsymbol{v}_{l-1} = \boldsymbol{R}_l \boldsymbol{v}_l, \qquad (15)$$
where $\boldsymbol{R}_l$ is the restriction operator. The restriction of the bilinear form on the other hand can be achieved in two commonly employed manners. One approach uses the restriction and prolongation operators as
$$\boldsymbol{A}_{l-1} = \boldsymbol{R}_l \boldsymbol{A}_l \boldsymbol{R}_l^T, \qquad (16)$$
which we refer to as `RAP` method. Another approach obtains the coarse matrices through the computation of the bilinear form using the function space of the coarse grid including a mesh-based estimation of the stabilization parameter. We refer to this approach as `assembly` method. The existence of the mesh-dependent stabilization parameter in the formulation in the previous sections sharply sets these approaches apart in the context of multilevel methods. This effect is investigated in Section 3.

**Smoothers** We use a multiplicative Schwarz-type smoother[31,32] for the treatment of cut cells on the finest mesh. The smoother consists in the solution of a series of local subdomain problems: we create a cell-based subdomain for every cell cut by the physical boundary that includes all degrees of freedom associated with the cell. A single-DoF subdomain is created for all other degrees of freedom which do not appear in a cell-based subdomain.[32] Given these $n_{sd}$ subdomains, the smoother operator can be written as
$$\boldsymbol{S} = (\boldsymbol{R}_{s,n_{sd}}^T \boldsymbol{A}_{s,n_{sd}}^{-1} \boldsymbol{R}_{s,n_{sd}}) \ldots (\boldsymbol{R}_{s,1}^T \boldsymbol{A}_{s,1}^{-1} \boldsymbol{R}_{s,1}) \quad \text{with} \quad \boldsymbol{A}_{s,i} = \boldsymbol{R}_{s,i} \boldsymbol{A} \boldsymbol{R}_{s,i}^T, \qquad (17)$$
where $\boldsymbol{R}_{s,i}$ is the restriction operator of the $i$-th subdomain and $\boldsymbol{A}_{s,i}$ is the local block of the $i$-th subdomain.

For the `assembly` approach on coarser meshes, the same smoother strategy could be applied. For the `RAP` approach, however, the formation of the coarse grid matrices is purely algebraic. Therefore, the employed smoothers usually are chosen to be algebraic methods as well. As the main target in our studies is a comparision of the approaches, we therefore employ a damped Jacobi smoother on coarse meshes in both cases. Please note that a Schwarz smoother on assembled coarse matrices likely provides even better convergence rates than the ones presented for comparison reasons in Section 3.

**Base solver** We use a direct LU solver on the base level. Given the high computational complexity of direct solvers[27], we note that it is important that the size of the base problem is kept relatively small, as discussed in Section 3.

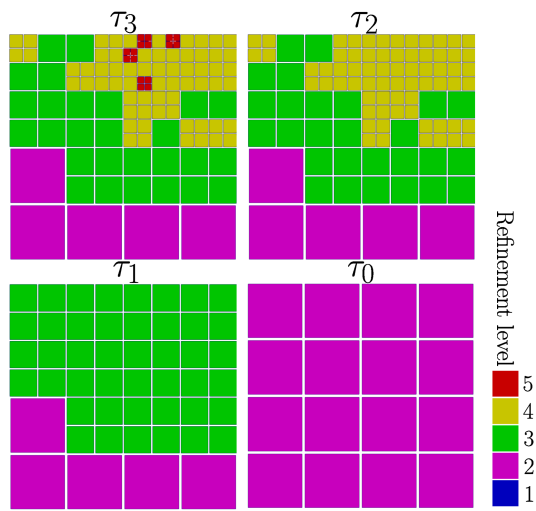

Figure 3: The hierarchy of grids constructed from a randomly-refined fine grid, $\tau_3$. Cells are colored based on their refinement level

Table 1: Mesh hierarchy: the coarse mesh $\tau_0$ is constructed by five levels of uniform refinement on the unit square. For $k \geq 1$, the meshes $\tau_k$ are constructed by adaptively refining its coarser parent $\tau_{k-1}$ one level towards the boundary of the physical domain. $n_{\texttt{DoF}}$ is the number of degrees of freedom in the mesh. $\lambda_g$ and $\lambda_l^{\text{mean}}$ are the global and mean local eigenvalues, respectively

| Mesh | $n_{\texttt{DoF}}$ | $\lambda_g$ | $\lambda_l^{\text{mean}}$ |
| --- | --- | --- | --- |
| $\tau_0$ | 1089 | - | - |
| $\tau_1$ | 1209 | 1279.3 | 312.23 |
| $\tau_2$ | 1657 | 13343.2 | 912.02 |
| $\tau_3$ | 2637 | 22770.2 | 2457.70 |
| $\tau_4$ | 4661 | 310834.0 | 4188.30 |
| $\tau_5$ | 8665 | 136533.0 | 17873.99 |
| $\tau_6$ | 16953 | - | 30835.11 |
| $\tau_7$ | 33269 | - | 62312.00 |
| $\tau_8$ | 65689 | - | 161481.18 |

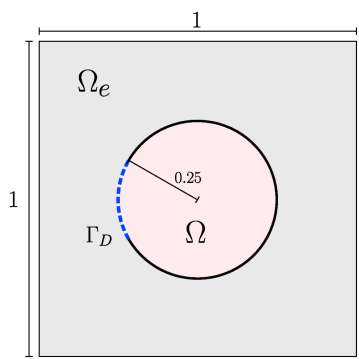

Figure 4: A circular physical domain is embedded in a unit square. A constant-valued inhomogeneous Dirichlet boundary is imposed on an arch to the left of the physical domain and the rest is left as homogeneous Neumann boundary

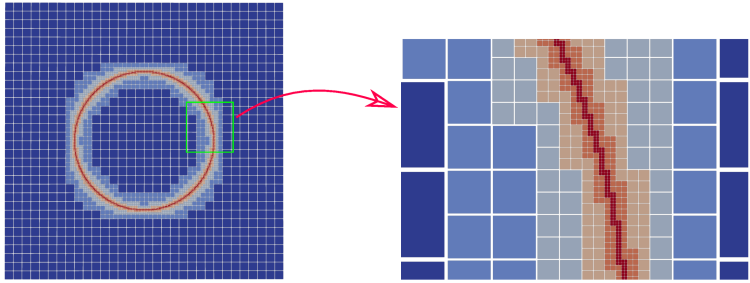

Figure 5: The discretization of the domain for $\tau_5$, i.e., five levels of uniform refinement on the unit square followed by five levels of adaptive refinement towards the boundary of the physical domain

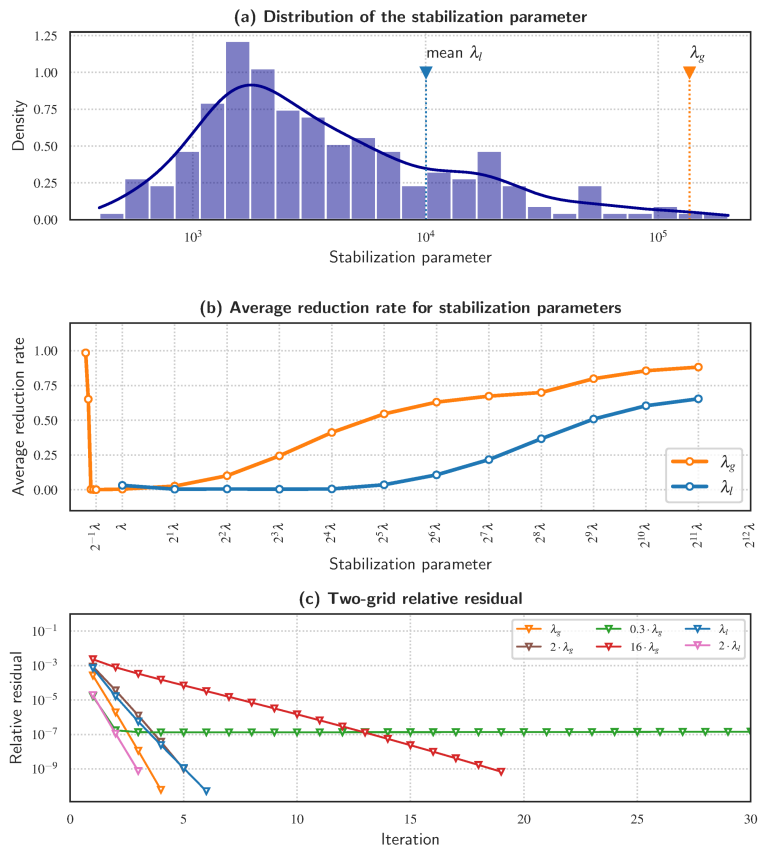

Figure 6: Two-grid solver on mesh $\tau_5$. (a) The distribution of the local stabilization parameter computed from the generalized eigenvalue problem and the global stabilization parameter, (b) the average reduction rate for variations in the global and local stabilization parameters and (c) the convergence of the iterative solver in terms of relative residual for selected cases. $\lambda_g$ and $\lambda_l$ denote the computed stabilization parameters using the global and local approaches, respectively. In (b), $\lambda$ denotes the element-wise constant stabilization parameter in the case of the local approach and the global stabilization parameter in the case of the global approach. Note that the solver does not converge for values smaller than the computed value in the local approach.

# 3 Numerical studies

For the numerical experiments, we consider a circular physical domain that is embedded inside a square as shown in Figure 4, where a Poisson problem with Dirichlet and Neumann boundary conditions is solved. An inhomogeneous Dirichlet boundary condition with a constant value of 1.0 is imposed on an arch section to the left of the domain and the rest is left as homogeneous Neumann boundary as demonstrated in Figure 4. Starting from the unit square, the refinement criteria includes five levels of uniform refinement followed by a number of adaptive refinement levels towards the boundary of the physical domain, leading to a hierarchy of nested grids that are used throughout this section. The grid hierarchy is detailed in Table 1. The discretization of the domain by mesh $\tau_5$ is shown in Figure 5. In addition, eight levels of adaptive quadrature are used for volume integration.

The multigrid solver detailed in Section 2.3 is used as a solver for the discretized problem. A V-cycle with three pre- and post-smoothing steps is employed. Relative residual with a threshold of $10^{-9}$ is taken as the convergence criterion. The penalization parameter $\alpha$ is $10^{-10}$. The stabilization parameter in Nitsche's formulation is computed through either the local, element-wise or the global approach detailed in Section 2.2. The computed global and local parameters from the generalized eigenvalue problems are designated with $\lambda_g$ and $\lambda_l$, respectively. An in-house code is used to implement the finite cell approach, the multigrid algorithm and to run the studies. p4est[34] and PETSc[36] are used for octree and linear algebra functionalities.

## 3.1 Stabilization parameter influence

In order to investigate the sensitivity of the iterative method to the stabilization parameter in each scheme, a number of variations from the computed parameter are examined: the global or the local stabilization parameter is multiplied by different factors to analyze the behavior for values smaller and larger than the ones suggested by the estimates obtained through the solution of the eigenvalue problems described in Section 2.2. The two-grid method on mesh $\tau_5$ is used in order to better isolate the effect of the stabilization parameter. The convergence of the solver as well as the distribution of the stabilization parameter are shown in Figure 6. Please note that the global stabilization parameter is not necessarily larger than the largest local parameter, as shown in Figure 6(a), i.e., the computed global stabilization parameter may be smaller than the local parameter for some of the cut cells.

For the global stabilization parameter scheme, the solution achieves rapid convergence in the vicinity of the computed value as can be seen in Figure 6(b). However, lower reduction rates can be obtained at smaller stabilization parameters than the one computed from the generalized eigenvalue problem. This trend continues up to the point that the stability of the solution is lost. At this point, we observe that the solution could either asymptotically reach a plateau or diverge. An example of such behavior is shown in Figure 6(c) for a stability parameter around three times smaller than the computed global value, where the solution reaches a plateau at approximately $10^{-7}$ and slowly diverges. Further reduction of the global stabilization parameter leads to divergence of the solver. It can be inferred that the most desirable behavior in the solution for this scheme is achieved for the smallest value of the global parameter that is large enough to ensure stability. An empirical choice of $2\lambda_g$ suggested by a number of studies[17,18] does not seem to have an obvious advantage in terms of the convergence rate of the solver for the global estimate.

For the local stabilization parameter scheme, the solution behavior is shown for variations of the stabilization parameter in Figure 6(b). It can be observed that the smallest reduction rates occur in the proximity of the computed local parameters. Although the solution quickly becomes unstable below the computed values, it remains consistent over a larger span compared to the global parameter. It can be inferred that the computed local stability parameter is close to the minimum value that is required for the coercivity of the problem. An empirical value of $2\lambda_l$ seems to be a reasonable choice in this case.

Although both methods achieve virtually identical convergence rates at the computed values, the average reduction rate achieved by the local scheme is considerably smaller than its global counterpart at higher relative parameters, and the iterative solver with the local scheme exhibits a better convergence behavior for the majority of the interval. Furthermore, the iterative method is more robust with regards to variations in the stabilization parameter in this scheme and the minimum reduction rate remains relatively constant over a rather large span of variation. The distribution of the local stabilization parameter in Figure 6(a) can shed some light on this observation. It can easily be seen that there is an enormous gap between the distribution of the stabilization parameter in the local scheme and the computed global value. Intuitively, the global value can be interpreted as an estimate for the minimum value that satisfies the coercivity condition on all elements, whereas the local values only need to satisfy the condition in a local domain. This freedom allows for a much more accurate distribution of the stability parameter that leads to better performance of the solver.

While other configurations of the solver, e.g., more pre/post smoothing steps certainly change the absolute value of measures such as the average reduction rate and similarly the number of required iterations, we did not observe fundamentally different behavior in the convergence of the solver for such variations.

## 3.2 Multigrid convergence

The implications of the mesh dependence of the stabilization parameter on a multigrid solver are analyzed through a mesh study in this section. The mesh hierarchy in Table 1 is employed. The mesh $\tau_0$ is used as the base problem for all tests, and is solved down to machine accuracy using a direct solver. Therefore, larger problems employ a deeper grid hierarchy. Local and global estimation of the stabilization parameter along with two approaches to the computation of the coarse grid matrices, namely `RAP` and `assembly`, lead to four possible cases. The configuration of the geometric multigrid solver is identical to the previous example on each level. The Schwarz smoother is applied to the finest grid and a damped Jacobi smoother is applied on all coarser levels. $2\lambda_g$ and $2\lambda_l$ are used in the global and local schemes, respectively. On the three finest problems, only the local scheme is applied as the computational cost of the global scheme becomes exceedingly inhibiting. The average reduction rate of the multigrid solver for each case is given in Figure 7(a).

It can be seen that for both methods of obtaining the system matrix, namely `RAP` and `assembly`, the local estimation method outperforms its global counterpart. The difference increases with levels of hierarchy. In each case of local and global estimation, the `assembly` methods achieves a smaller reduction rate than its `RAP` counterpart. A significant gap in the convergence rate can be observed between the combinations at higher levels of the hierarchy.

When the coarse grid matrices are assembled, the convergence rate achieved by the global estimation scheme remains relatively close to the local estimation scheme up to the first few levels, after which the convergence rates start to deviate. The underlying mechanism of this behavior is twofold. In this case, the global stabilization parameter is estimated on each level of the hierarchy and is therefore suitable for that discretization. This can be inferred as the main reason why the global scheme performs relatively well up to a certain level in this configuration. The deterioration of the convergence rate can be imputed to slightly inferior convergence rates at each level that accumulate in the multilevel configuration. This is only exacerbated by the enlargement of the hierarchy depth.

An analogous behavior is observed in the case where coarse grid matrices are obtained through restriction. A more accurate distribution of the stabilization parameter through the local scheme on the fine grid leads to lower reduction rates, especially at higher levels of hierarchy.

In either estimation method of the stabilization parameter, smaller and more scalable convergence rates are achieved when using the `assembly` method for coarse matrices. This can be attributed to the fact that a more accurate estimation of the stabilization parameter at each level can achieve a better convergence rate. In the `RAP` method, this estimation is performed only on the finest grid, and might not be reasonable for coarser grids. This can be more easily perceived by looking at the distribution of the stabilization parameters on different depths of the hierarchy. The distribution of the stabilization parameter is shown in Figure 7(b) for two sample grids. It can clearly be seen that the stabilization parameter has a highly grid-dependent distribution. If the distribution from one grid were to be used on the other, it would either lead to extremely slow convergence or divergence according to the results in the previous section. Since the majority of cells are larger on coarser grids, the stabilization estimates from finer grids are normally sufficiently large for the coarse grid. The convergence rates and relative robustness of the `RAP` method may be justified by this mechanism. However, we would like to bring a possible scenario into attention. The stabilization estimate depends not only on the grid size, but also on the cut configuration. Specifically, unfortunate cut elements, where the physical boundary intersects a small portion of the cell lead to large stabilization parameters. It is thus possible that the required stabilization parameter on a coarser grid be larger than on the fine grid. This is manifested for example between $\tau_4$ and $\tau_5$ in Table 1, where the coarser grid, $\tau_4$, requires a larger stabilization parameter than the finer grid, $\tau_5$, in the global estimate although on average, the required stabilization parameter is smaller for $\tau_4$ in the local estimate. Furthermore, Figure 8 illustrates a possible scenario in which the cut configuration is deteriorated on the coarse grid although the grid size is larger. In conclusion, it is observed that a dedicated estimate on every mesh level through the `assembly` method can lead to significant improvement in the convergence rate of the solver.

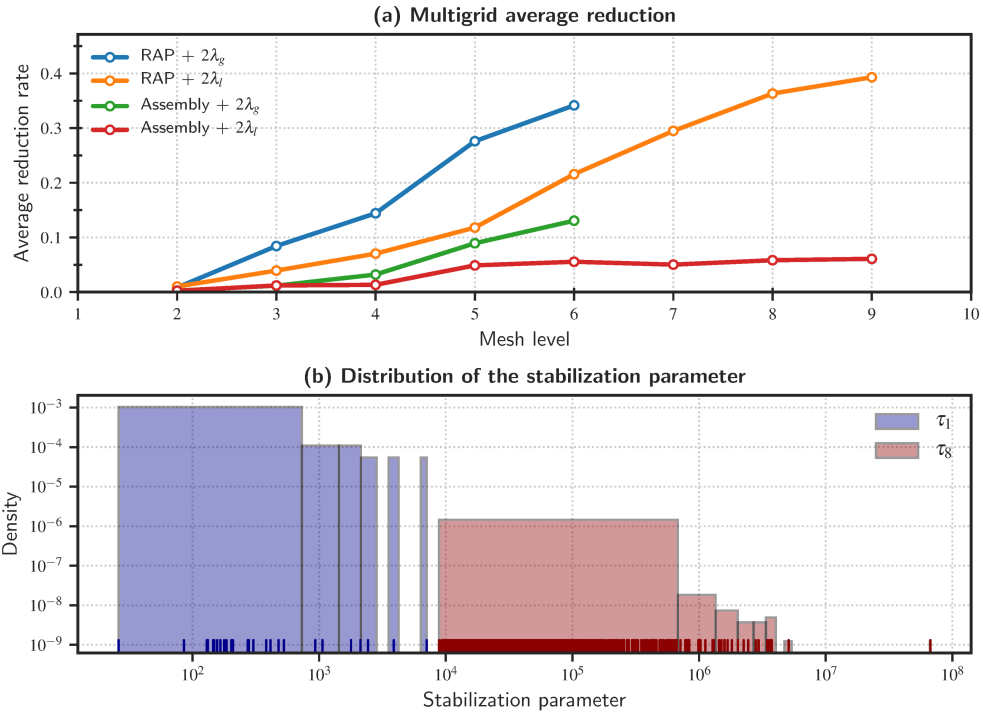

Figure 7: Study for multigrid hierarchy: (a) Average reduction rate for the geometric multigrid solver for global and local estimation of the stabilization parameter and two approaches of obtaining the coarse grid matrices, (b) the distribution of the local stabilization parameter on a coarse ($\tau_1$) and a fine ($\tau_8$) mesh

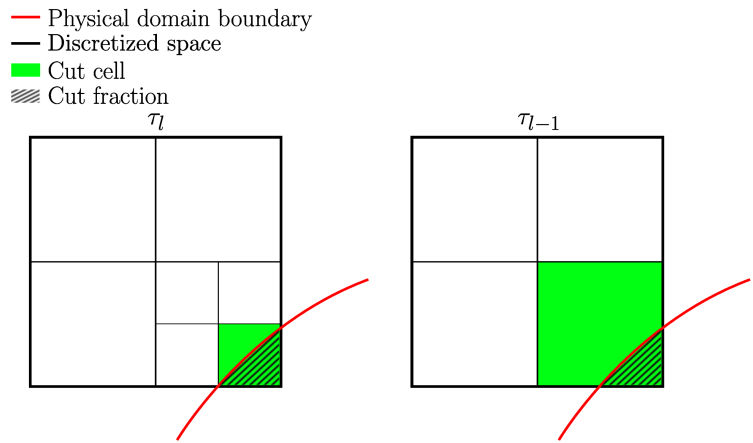

Figure 8: A possible scenario where the smaller cut fraction on the coarser grid could potentially lead to a larger stabilization parameter

# 4 Conclusions

We investigate the mesh-dependent stabilization parameter in Nitsche's method in the context of the geometric multigrid solution to immersed finite element formulations. We find that the stabilization parameter not only carries importance for ensuring the coercivity of the bilinear operator and thus the stability of the solution, but also significantly affects the performance of the iterative solver.

The stabilization parameter can in general be incorporated in the integration process either as a domain-wide or an element-wise constant. A good estimate for each approach can be formulated through generalized eigenvalue problems, leading to the global and local estimates, respectively. The global estimate obtains a constant parameter from a single generalized eigenvalue problem, whereas a series of smaller generalized eigenvalue problems lead to a distribution of the stabilization parameter in the local estimate. For a multigrid setting, the stabilization parameters can be computed separately for each mesh level if the coarse grid matrices are assembled. Using the `RAP` approach for coarse grid matrices, the stabilization parameter estimate is only computed on the finest mesh and then implicitly contained in coarser matrices via restriction.

For the stabilization parameter methods, we find that deviation from the values obtained from the generalized eigenvalue problems significantly affects the convergence rate of the solver. Specifically, the best behavior in the global estimate is achieved for the smallest global parameter that ensures stability of the solution, which can be several times smaller than the value computed from the generalized eigenvalue problem. For the local estimation, the solver tends to achieve its smallest reduction rates for a value of $2\lambda_l$-$4\lambda_l$. The solution quickly becomes unstable below the calculated distribution in this approach. However, it is more robust for larger relative values compared to the global estimate. The reduction rate of the solver deteriorates for overestimated parameters in both approaches. We find that the solver tends to achieve better convergence rates with the local estimate method, especially for deeper hierarchies.

The discretization dependence of the stabilization parameter influences the performance of the multilevel solver. Forming coarse grid matrices through the evaluation of the bilinear form on the coarse function space, an individual stabilization parameter estimate on each level of the hierarchy is obtained which we find to provide better convergence rates. For the `RAP` approach, the finest level estimate of the stabilization parameter is implicitly used also on coarser meshes which deteriorates the convergence rates, in particular for deeper mesh hierarchies.

We achieve the fastest and most scalable convergence rates for a local stabilization parameter estimate in conjunction with assembled matrices on each level of the hierarchy. The local estimate seems to be a favorable approach to obtaining the stabilization parameter, considering the increasing costs for generalized eigensolvers for the global estimate and the better performance achieved by the local estimate. Computing a separate stabilization parameter on every level of the multilevel hierarchy via the evaluation of the bilinear form on coarse spaces seems favorable as it produces level-adjusted stabilization parameters in the coarse matrices which result in better convergence rates for the multigrid solver in our studies.

# Acknowledgments

Financial support was provided by the German Research Foundation (Deutsche Forschungsgemeinschaft, DFG) in the framework of subproject C4 of the collaborative research center SFB 837 Interaction Modeling in Mechanized Tunneling, grant number 77309832. This support is gratefully acknowledged. We also gratefully acknowledge the computing time on the computing cluster of the SFB 837 and the Department of Civil and Environmental Engineering at Ruhr University Bochum, which has been employed for the presented studies.